\documentclass{IET}

\usepackage{xcolor}
\usepackage[colorlinks, 
linkcolor={red!50!black}, 
citecolor={blue!50!black}, 
urlcolor={blue!80!black}]{hyperref}

\usepackage[numbers,sort&compress]{natbib}

\usepackage{algorithm}
\usepackage{algpseudocode}
\usepackage{calc}

\newcommand{\A}{\ensuremath{\mathcal{A}}}

\newcommand{\range}{\ensuremath{R}}

\newcommand{\scanpattern}{\ensuremath{\mathcal{S}}}

\newcommand{\grid}{\ensuremath{G}}
\newcommand{\cover}{\ensuremath{C}}
\newcommand{\covers}{\ensuremath{\mathcal{C}}}

\newcommand{\bc}{\ensuremath{\bold{c}}}
\newcommand{\bC}{\ensuremath{\bold{C}}}
\newcommand{\bA}{\ensuremath{\bold{A}}}
\newcommand{\bT}{\ensuremath{\bold{T}}}
\newcommand{\bx}{\ensuremath{\bold{x}}}

\newcommand{\node}{\ensuremath{N}}
\newcommand{\nodes}{\ensuremath{\mathcal{N}}}

\DeclareMathOperator*{\argmin}{argmin}
\DeclareMathOperator*{\com}{com}

\begin{document}

\title{Branch-and-Bound Method for Just-in-Time Optimization of Radar Search Patterns}

\author[1,2,*]{Yann Briheche}
\author[1]{Frederic Barbaresco}
\author[2]{Fouad Bennis}
\author[2]{Damien Chablat}

\affil[1]{THALES AIR SYSTEMS, Voie Pierre-Gilles de Gennes, 91470 Limours, France}
\affil[2]{Laboratoire des Sciences du Num{\'e}rique de Nantes, UMR CNRS 6004, 44321 Nantes, France}
\affil[*]{\href{mailto:corresponding.author@second.com}{yanis.briheche@thalesgroup.com}}



\maketitle

\section{Introduction \& Context}\label{sec1}

Set covering is a well-known problem in combinatorial optimization. The objective is to cover a set of elements, called the universe, using a minimum number of available covers. The theoretical problem is known to be generally NP-difficult to solve \cite{Vazirani2001}, and is often encountered in industrial processes and real-life problem. In particular, the mathematical formulation of the set cover problem is well-suited for radar search pattern optimization of modern radar systems.

\begin{figure}
	\centering
	\includegraphics[width=0.70\linewidth]{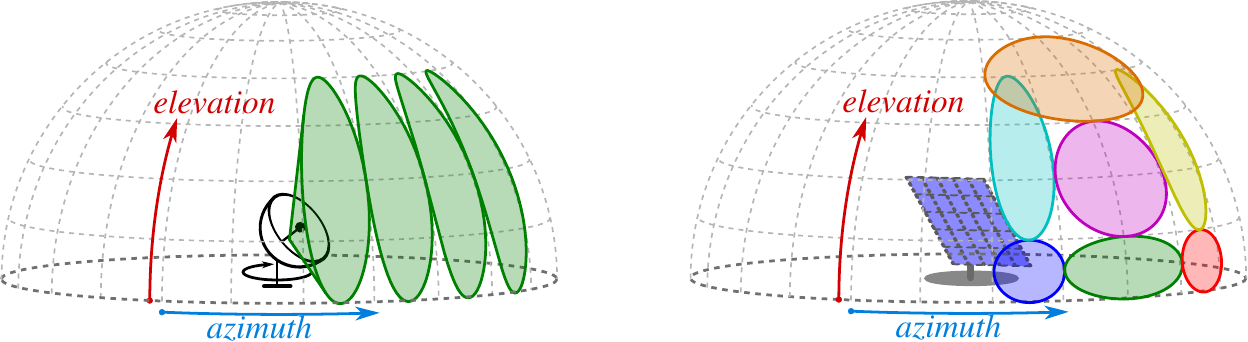}
	\caption{Radar search pattern for a rotating radar (left) and a modern fixed-panel radar (right)}
	\label{modernradars}
\end{figure}

Electronic scanning and numerical processing allow modern radars to dynamically use bi-dimensional beam-forming, giving them great control on the radar search pattern. While traditional rotating radars search the space sequentially along the azimuth axis and reproduce at each azimuth the same pattern along the elevation axis, modern radars can optimize the search pattern along both axis simultaneously (Fig. \ref{modernradars}). Those new possibilities requires a more sophisticated formulation for the optimization problem of the radar search pattern. Approximation of radar search pattern optimization as a set cover problem offers a flexible yet powerful formulation \cite{Briheche2016}, capable of accounting radar specific constraints (localized clutter, adaptive scan-rate updates, multiple missions) without changing the underlying mathematical structure of the optimization problem.

Various optimization algorithms have been proposed for solving the set cover problem: exact methods such as branch-and-bound combined with relaxations methods, and approximation algorithms such as greedy algorithm, simulated annealing and genetic algorithms (see \cite{Yelbay2015} for a recent survey of those methods). In practice, problems of reasonable size can be efficiently solved by branch-and-bound exploration using linear relaxation for lower bound estimation. More importantly, branch-and-bound features interesting characteristics, making it particularly fit for producing just-in-time solutions, for example for radars in operational situation.

\section{Problem Statement}

\subsection{Definition}

Let $\grid=\{g_{m,n}\}$ be the set representation of a finite bi-dimensional $M$-by-$N$ regular grid with 
\begin{itemize}
\item each element $g_{m,n}$ representing a cell indexed by $(m,n)\in [0,M[ \times [0,N[ \ \subset \mathbb{N}^2$.
The grid contains $MN$ cells.
\item each couple $(m,n)$ representing the node at the intersection of the $m$-th horizontal line and the $n$-th vertical line with $(m,n)\in [0,M] \times [0,N] \ \subset \mathbb{N}^2$. The grid has $(M+1)(N+1)$ nodes.
\end{itemize}

On the grid, a rectangular cover is a subset of elements included in a rectangle, uniquely defined by its upper left corner node $(m_0,n_0)$ and its lower right corner node $(m_1,n_1)$, such that $0\leq m_0 < m_1 \leq M$ and $0 \leq n_0 < n_1 \leq N$. The set representation of a cover defined by corners $(m_0,n_0)$ and $(m_1,n_1)$ is:
$$\cover=\{g_{m,n},\ (m,n) \in [m_0, m_1[ \times [n_0, n_1[\}$$
For example, in cover $\cover_7$ (Fig. \ref{setcoverexample}), the corners are $(m_0,n_0)=(0,1)$ and $(m_1,n_1)=(1,2)$.\\

\noindent Let $\covers=\{\cover_1,\dots,\cover_D\}$ be a collection of $D$ rectangular covers on $\grid$.\\
Let $T_\cover \in \mathbb{R^+}$ be the associated cost of cover $\cover \in \covers$ (also noted $T_i$ for cover $\cover_i$).\\
Find a minimum cost sub-collection $\scanpattern \subset \mathcal{C}$ covering all cells of grid $\grid$.\\

\subsection{Example}

\begin{figure}
	\centering
	\includegraphics[width=\linewidth]{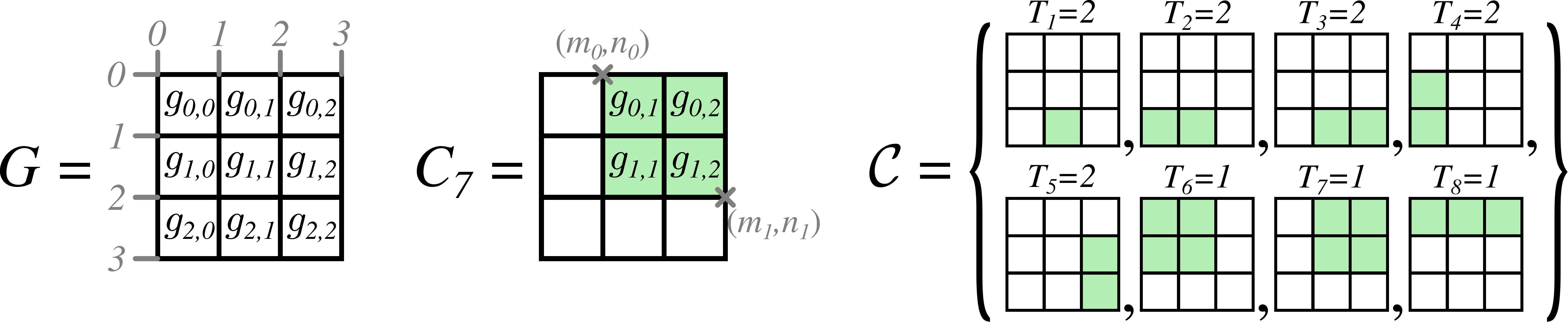}
	\caption{Grid $\grid$ to cover (left), a cover (center) and the collection $\covers$ of available covers (right)}
	\label{setcoverexample}
\end{figure}

There are six available covers such as in Figure \ref{setcoverexample} to cover $\grid$: 
\begin{itemize}
	\item $\scanpattern_1=\{\cover_1,\cover_4,\cover_5,\cover_6,\cover_7\}$ is a valid sub-optimal covering collection with total cost $T_1+T_4+T_5+T_6+T_7=8$.
	\item $\scanpattern_2=\{\cover_2,\cover_3,\cover_6,\cover_7\}$ is a valid optimal covering collection with total cost $T_2+T_3+T_6+T_7=6$, as there are no solution with total cost $5$ or less.
	\item $\scanpattern_3=\{\cover_2,\cover_5,\cover_6,\cover_8\}$ is another optimal covering collection, thus an optimal solution is not necessarily unique.
\end{itemize}

The optimization formulation of this set cover problem can be written as:

\begin{equation}
\begin{array}{rl}
\min & \sum_{\cover \in \scanpattern} T_\cover \\
s.t. & \forall g_{m,n} \in \grid, \exists \cover \in \scanpattern,\ g_{m,n} \in  \cover \\
	 & \scanpattern \subset \covers
\end{array}\label{setcoverequation}
\end{equation}

\begin{figure}
	\centering
	\includegraphics[width=0.80\linewidth]{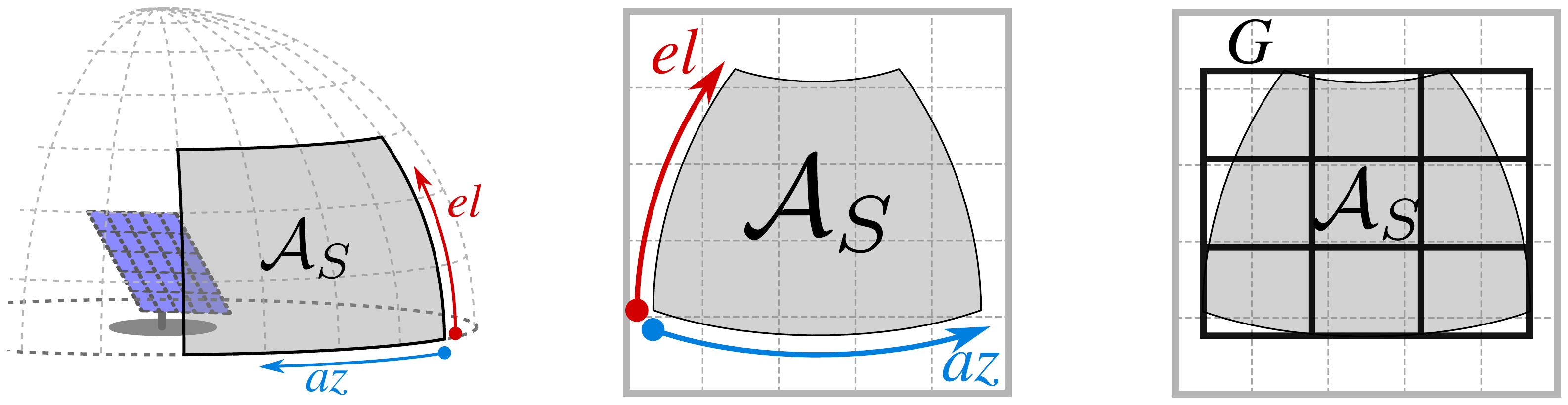}
	\caption{Surveillance area $\A_S$ (left), its projection in direction cosines (center) and the surveillance grid (right)}
	\label{surveillancearea}
\end{figure}

\begin{figure}
	\centering
	\includegraphics[width=0.80\linewidth]{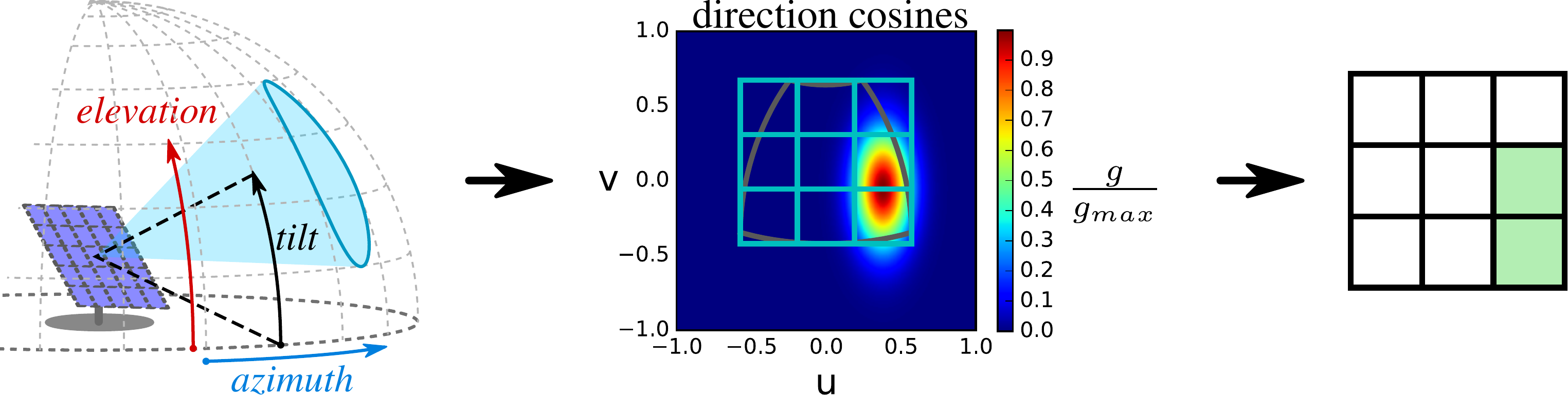}
	\caption{Radar detection beam (left), its radiation pattern (center) and the associated cover (right)}
	\label{beam}
\end{figure}

In the case of radar search pattern application, the grid $\grid$ represents the surveillance area (Fig. \ref{surveillancearea}), while each cover $\cover \in \covers$ represents a radar beam and its detection area on the grid $\grid$ (Fig. \ref{beam}). The associated cost $T_\cover$ is the duration required to emit the radar signal, then receive and process the echo. The total cost of collection of radar beams is the time required to emit all beams in sequential order, as the radar cannot emit simultaneously several beams.

\subsection{Combinatorial complexity}

A rectangle cover is uniquely define by its upper left and lower right corners. Those corners are mathematically defined by choosing two values $m_0$ and $m_1$ among the $M+1$ horizontal lines, and two values $n_0$ and $n_1$ among the $N+1$ vertical lines on the grid. Thus there are at most $$\binom{M+1}{2}\binom{N+1}{2}=\frac{MN(M+1)(N+1)}{4} = O(M^2N^2) $$ possible distinct rectangles on a $M$-by-$N$ grid. And so the maximum number of possible sub-collections of rectangular covers on the grid is $2^{MN(M+1)(N+1)/4}$.

Even for a $10$-by-$10$ grid, which is relatively small, the number of possible sub-collections is approximately $10^{900}$, which is far too big to allow the use of brute-force exploration.

\section{Integer Programming}

\subsection{Problem Formulation}

The set cover problem can be written as an integer program by using matrix formulations. We represent each cover $\cover \in \covers$ as a binary $M$-by-$N$ matrix noted $\bold{C}$, or as a binary vector of length $MN$ noted $\bold{c}$:
$$
\bold{C}(m,n)=
\bold{c}(m+Mn)=
\begin{cases}
1 \text{ if } g_{m,n} \in \cover\\  
0 \text{ otherwise}
\end{cases} 
$$

\begin{figure}
	\centering
	\includegraphics[width=0.90\linewidth]{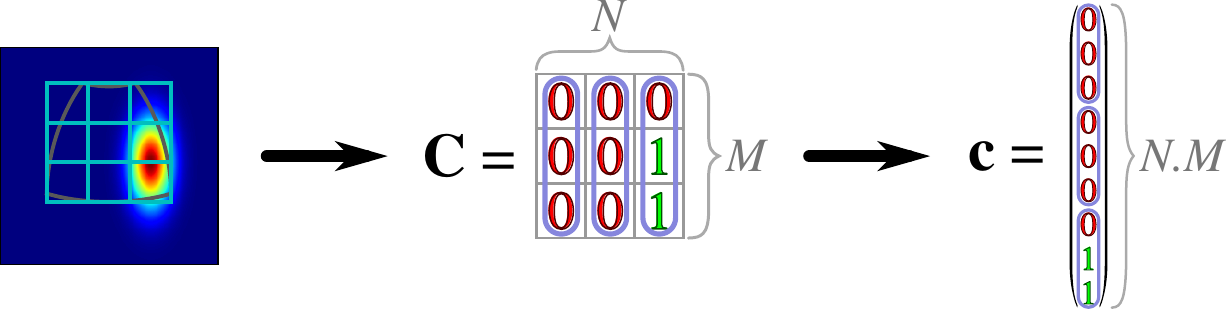}
	\caption{Radar detection beam (left), its binary matrix representation (center) and its binary vector representation (right)}
	\label{covervector}
\end{figure}

For each cover $\cover_i \in \covers$, let $x_i \in \{0,1\}$ be the binary selection variable of cover $\cover_i$, such that the vector $\bold{x}=(x_1,\dots,x_D)\in \{0,1\}^D$ represents the sub-collection $\scanpattern=\{\cover_i \in \covers \ s.t.\ x_i=1\}$, containing the chosen covers.

Let $\bold{T}=(T_1\cdots T_D)^T$ be the cost vector and let $$\bold{A}=
\begin{pmatrix}
\bold{c}_1 & \cdots & \bold{c}_D
\end{pmatrix}
=
\begin{pmatrix}
\bold{C}_1(0,0) & \cdots & \bold{C}_D(0,0) \\
\bold{C}_1(1,0) & \cdots & \bold{C}_D(1,0) \\
\vdots & \ddots & \vdots\\
\bold{C}_1(m,n) & \cdots & \bold{C}_D(m,n) \\
\vdots & \vdots & \vdots\\
\end{pmatrix}$$
be the cover matrix.

Then the set cover problem (\ref{setcoverequation}) can be written as the following integer program:
\begin{equation}
\begin{array}{rl}
	\min & \bold{T}^T.\bold{x}                  \\
	s.t. & \bold{A}\cdot \bold{x} \geq \bold{1} \\
	     & \bold{x} \in \{0,1\}^D
\end{array}\label{integerprogram}
\end{equation}
where $\bold{1}$ is the vector $(1 \cdots 1)$ of length $MN$. As an example, the set cover problem represented in Figure \ref{setcoverexample} can described by the following Equation:
\begin{equation}
\bold{A}=
\begin{pmatrix}
	0 & 0 & 0 & 0 & 0 & 1 & 0 & 1 \\
	0 & 0 & 0 & 1 & 0 & 1 & 0 & 0 \\
	0 & 1 & 0 & 1 & 0 & 0 & 0 & 0 \\
	0 & 0 & 0 & 0 & 0 & 1 & 1 & 1 \\
	0 & 0 & 0 & 0 & 0 & 1 & 1 & 0 \\
	1 & 1 & 1 & 0 & 0 & 0 & 0 & 0 \\
	0 & 0 & 0 & 0 & 0 & 0 & 1 & 1 \\
	0 & 0 & 0 & 0 & 1 & 0 & 1 & 0 \\
	0 & 0 & 1 & 0 & 1 & 0 & 0 & 0
\end{pmatrix}
\text{, }
\bold{T}=
\begin{pmatrix}
2 & 2 & 2 & 2 & 2 & 1 & 1 & 1
\end{pmatrix}^T
\text{ and }
\bx=
\begin{pmatrix}
x_1 \\ x_2 \\ x_3 \\ x_4 \\ x_5 \\ x_6 \\ x_7 \\ x_8
\end{pmatrix}
\label{setcoverexampleequation}
\end{equation}

In their general form, integer programs are NP-hard to solve. Intuitively, this means that solving those problems is difficult, and requires some form of exhaustive enumeration of all possible solutions, whose number is often exponential in respect to the problem size.

\subsection{Linear Relaxation}

The linear relaxation of an integer program can be obtained by relaxing the integrality constraint of (\ref{integerprogram}) into a positivity constraint, allowing the variables $(x_i)_{1\leq i \leq D}$ to take continuous values in $[0,1]$:
\begin{equation}
\begin{array}{rl}
\min & \bold{T}^T.\bold{x}                  \\
s.t. & \bold{A}\cdot \bold{x} \geq \bold{1} \\
& \bold{0} \leq \bold{x} \leq \bold{1}
\end{array}\label{linearprogram}
\end{equation}
Any valid solution of the integer program is also a valid solution of its linear relaxation. Consequently the optimal value of the linear relaxation is inferior to the optimal value of the integer program, since an optimal solution of the integer program is a valid solution of the linear relaxation.

Note that the constraint $\bold{x}\leq \bold{1}$ is in fact unnecessary, since the problem
\begin{equation}
\begin{array}{rl}
\min & \bold{T}^T.\bold{x}                  \\
s.t. & \bold{A}\cdot \bold{x} \geq \bold{1} \\
& \bold{0} \leq \bold{x}
\end{array}\label{linearprogramreduced}
\end{equation}
has the same optimal solutions as (\ref{linearprogram}). Intuitively, in the linear relaxation, a cell is going to be covered by a sum of ``fractional" covers (with $x_i<1$), or as at least one integer cover (with $x_i=1$) and thus has no need for covers with $x_i>1$. 

We formalize mathematically this idea in the following proof:
\begin{itemize}
	\item Let $\bold{x}^b=(x_1^b \cdots x_D^b)$ be an optimal solution of (\ref{linearprogramreduced}). Let $\bold{y}=(y_1 \cdots y_D) \text{ with }y_i=\min\{x_i^b,1\}$.\\
	Immediately we have $\bT^T \bold{y} \leq \bT^T \bx^b$. Since $\bold{x}^b$ is a valid solution of (\ref{linearprogramreduced}): $$\forall (m,n),\sum_i x_i^b \cover_i(m,n) = \sum_{x_i^b\leq 1} x_i^b \cover_i(m,n) + \sum_{x_i^b> 1} x_i^b \cover_i(m,n) \geq 1$$
	\begin{itemize}
		\item[-] Case 1: $\sum_{x_i^b>1} x_i^b \cover_i(m,n)=0$
		$$\sum_i y_i \cover_i(m,n) \geq \sum_{x_i^b\leq 1} x_i^b \cover_i(m,n) = \sum_{x_i^b\leq 1} x_i^b \cover_i(m,n) + \sum_{x_i^b> 1} x_i^b \cover_i(m,n) \geq 1$$
		\item[-] Case 2: $\sum_{x_i^b>1} x_i^b \cover_i(m,n)>0$\\
		$\exists j \text{ s.t. }\ x_j^b>1 \text{ and } \cover_j(m,n)>0$, and since $\cover_j(m,n) \in \{0,1\}$, $\cover_j(m,n)=1$, thus:
		$$ \sum_i y_i \cover_i(m,n) \geq y_j \cover_j(m,n) \geq \cover_j(m,n) \geq 1 $$
	\end{itemize}
	So $\bold{y}$ is a valid solution of (\ref{linearprogramreduced}), and since $\bT^T \bold{y} \leq \bT^T \bx^b $, $\bold{y}$ is also an optimal solution of (\ref{linearprogramreduced}).\\
	From this, we deduce $\bT^T \bold{y} = \bT^T \bx^b $ and with $\bT > \bold{0}$, we deduce $\bx^b=\bold{y}$, so $\bx^b$ is a valid solution for (\ref{linearprogram}).
	\item Let $\bx^a$ be an optimal solution of (\ref{linearprogram}), $\bx^a$ is also a valid solution of (\ref{linearprogramreduced}), so $\bT^T \bx^b \leq \bT^T \bx^a$. We just showed that if $\bx^b$ is an optimal solution of (\ref{linearprogramreduced}) then it is also a valid solution of (\ref{linearprogram}), so we also have $\bT^T \bx^a \leq \bT^T \bx^b$, and thus we can conclude $\bT^T \bx^b = \bT^T \bx^a$.

	So any optimal solution of (\ref{linearprogram}) is an optimal solution of (\ref{linearprogramreduced}) and reciprocally, thus both problems have the same set of optimal solutions.
		
\end{itemize}

Furthermore, the positivity constraints $\bold{0} \leq \bold{x}$ can be integrated in the matrix formulation with
$$
\bold{R} = 
\begin{pmatrix} 
\bA \\ 
\bold{I} 
\end{pmatrix}
\text{ and }
\bold{d}=
\begin{pmatrix} 
\bold{1} \\ 
\bold{0}
\end{pmatrix}
$$
by rewriting the linear program as
\begin{equation}
\begin{array}{rl}
\min & \bold{T}^T.\bold{x}                  \\
s.t. & \bold{R}\cdot \bold{x} \geq \bold{d}
\end{array}\label{linearprogramcomplete}
\end{equation}
And the three formulations of the linear relaxation (\ref{linearprogram}), (\ref{linearprogramreduced}) and (\ref{linearprogramcomplete}) are equivalent.
%

The integer program representing our set cover problem and its linear relaxation have two more interesting properties:
\begin{itemize}
	\item Easily-checked feasibility: an integer program is feasible if there is at least one solution validating all constraints. It is possible that no valid solution exists if some constraints are conflicting, or if one constraint is impossible. In our case, feasibility is easy to check: the integer program as well as its linear relaxation are feasible if and only if $\bx_F=(1\cdots 1)$ is a feasible solution, i.e. $\bA \cdot \bx_F = \sum_{i=1}^{D}\bc_i \geq \bold{1}$:
	\begin{itemize}
		\item if $\bx_F$ is a valid solution, then the problem is feasible by definition.
		\item if $\bx_F$ is an invalid solution, then there is an invalidated constraint for $\bx_F$, i.e.:\\
		$\exists (m,n)$ s.t. $\sum_{i=1}^{D}\bC_i(m,n)<1$, and since $\forall(i,m,n)$, $\bC_i(m,n) \in \{0,1\}$,\\ $\exists (m,n)$ s.t. $\forall i ,\bC_i(m,n)=0 \Rightarrow \exists (m,n)$ s.t. $\forall (x_i)_{1\leq i \leq D}, \sum_{i=1}^{D}x_i\bC_i(m,n)=0<1$\\		
		In other words, $A$ has its $(m+Mn)$-th row filled with zeros, corresponding to a constraint which can be satisfied by no solution. 
		
		Intuitively, $\bx_F$ represents $\covers$, the collection of all available covers itself, and if it is an invalid solution, then there is a cell which cannot be covered. This can happen in a real system if there is a cell which cannot be scanned, because of an obstacle or because the radar has not enough power to achieve the desired detection range.
		
	\end{itemize}
	
	\item Boundedness: a recurring question for linear programs is whether they are bounded, that is whether the cost function is bounded (below for minimization) for valid solutions. In our case the cost function is positive and thus bounded below by $0$.
\end{itemize}

\subsection{Linear Programming}

\begin{figure}
	\centering
	\includegraphics[width=0.7\linewidth]{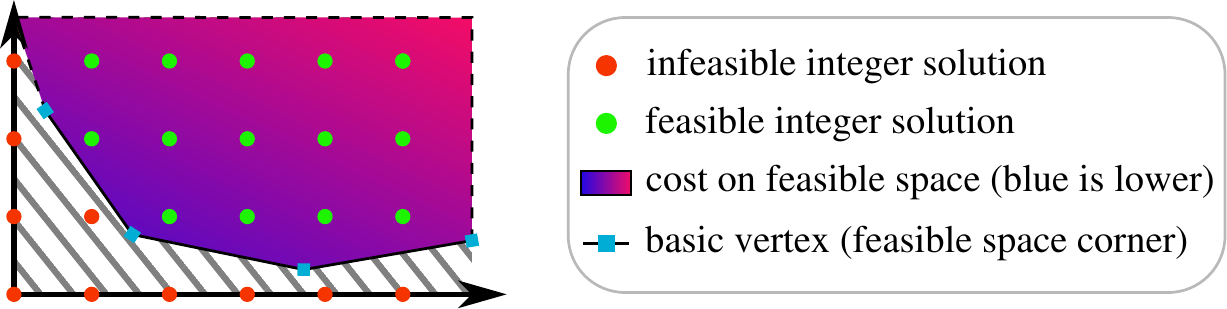}
	\caption{Decision space for 2D linear and integer programs}
	\label{lp_space}
\end{figure}

There are three important geometrical aspects describing the decision space of the integer and linear programs (Fig. \ref{lp_space}):
\begin{itemize}
	\item $\bT$ is the cost function gradient. The cost function is linear and its gradient is constant. $-\bT$ is the direction of maximum decrease of the cost function.
	\item $\bA$ is the cover matrix. Each row of $\bA$ correspond to a detection constraint on a cell of $\grid$. In the decision space, each constraint corresponds to an hyperplane, the limit between the halfspace of solutions validating the constraint and the halfspace of solutions violating the constraint. The intersection of those halfspace forms a convex polyhedron.
	\item The positivity constraint of the linear relaxation $ \bold{0} \leq \bold{x} \leq \bold{1}$ bounds the values of the valid solutions in the hypercube $[0,1]^D$.
	
	For the integer program, the integrality constraint $\bx \in \{0,1\}^D$ further reduces the set of valid solutions to the vertices of the hypercube $[0,1]^D$.
\end{itemize}


So the set of valid solutions for the linear relaxation is the intersection of the valid halfspaces for all constraints, and the hypercube $[0,1]^D$. Geometrically, it is a bounded convex polyhedron in $\mathbb{R}^D$, and can be described by its vertices (``corners"). Each vertex of this polyhedron is a point where at least $D$ hyperfaces of the polyhedron intersect, in other words, a point where $D$ constraints are tight. 

Such a point is called a basic solution (or basic vertex) of the linear program. It has been proved that if a linear program is bounded and feasible, then it has a basic optimal solution \cite{Matouek2006}. 

\begin{figure}
	\centering
	\includegraphics[width=0.7\linewidth]{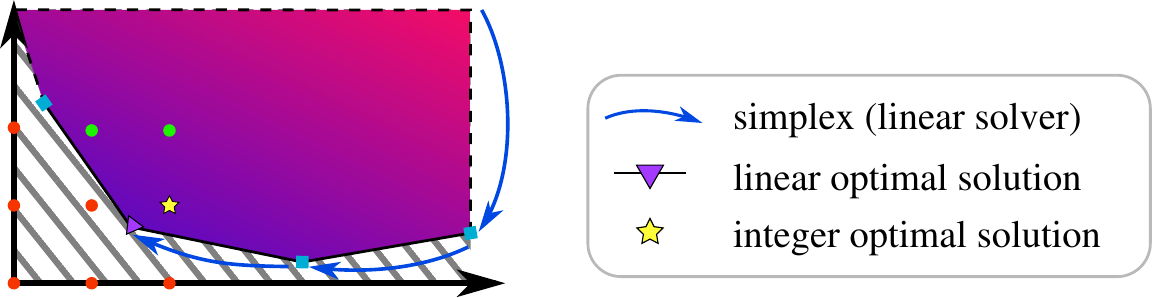}
	\caption{Illustration of Dantzig's simplex method for solving linear programs}
	\label{lp_solv}
\end{figure}

For a linear program, the polyhedron convexity allows the use of descent methods, such as Dantzig's simplex method, represented in Figure \ref{lp_solv}, which moves from vertex to vertex on the feasible polyhedron until it reaches an basic optimal solution, i.e. a vertex with no decreasing neighbor. However, this type of method generally cannot be used to solve integer programs, for which solutions are isolated points.

So for a given basic optimal solution, we have $D$ tight constraints. Let $B \leq MN$ be the number of tight detection constraints. If $D > B$ then we have $Z = D-B$ tight bound constraints. By considering formulation (\ref{linearprogramreduced}), we know that those $Z$ tight bound constraints are of the form $x_i \leq  0$, thus $x_i =  0$. The corresponding $Z$ variables are called \textit{non-basic} variables and are zeroes. The other $D-Z=B$ variables are called \textit{basic} variables and can be non-zero values. We reorganize the variables as $\bx^T=(\bx_B^T\ \bx_Z^T)$, where $\bx_B$ are the basic variables, and $\bx_Z$ are the non-basic variables. Thus we have $B$ tight detection constraints in $\bA$, such that 
$$\bA_B \bx_B = \bold{1}$$ 
where $\bA_B$ is the square $B$-by-$B$ submatrix of $\bA$ linking the basic variables $\bx_B$ to the tight detection constraints. Furthermore, $\bA_B$ is necessarily non-singular: since the hyperplanes of all constraints intersect, the constraints are linearly independent.

\subsection{Integral Program and Total Unimodularity}

So for any basic optimal solution of the linear program, there is a (possibly non-unique) square non-singular submatrix $\bA_B$ of $\bA$ such that $\bA_B \bx_B = \bold{1}$ (note that the reverse is not true, the condition is necessary, but not sufficient). Thus, $\bx_B=\bA_B^{-1} \bold{1}$ and $\bx^T=(\bx_B^T\ \bx_Z^T)=(\bx_B^T\ \bold{0})$.

So if $\bA_B^{-1}$ is an integral matrix (i.e. contains only integral values), then $\bx$ is also integral. This means that the linear program and the integer program share an optimal solution. The integrality of an invert matrix is determined by the determinant of its forward matrix: 
$$\det(\bA_B) \in \{-1,+1\} \Rightarrow \bA_B^{-1} = \frac{\com(\bA_B)^T}{\det(\bA_B)^{-1}}\text{ is integral}$$
A matrix $\bA$ is said to be \textit{unimodular} if $\det(\bA) \in \{-1,0,+1\}$. A matrix is said to be \textit{totally unimodular} if all its square submatrices are unimodular. An integer program whose constraint matrix is totally unimodular can be directly solved by linear relaxation, because the integer program and its linear relaxation have the same basic optimal solutions. In this case, the problem and its associated convex polyhedron are said to be integral. Geometrically, this means that all vertices of the polyhedron are integral points.

Total unimodularity is an important concept in combinatorial optimization, because it reduces integer programming to linear programming, which is theoretically an ``easier" problem. Linear programming is solvable in polynomial time, and very efficient practical algorithms exist.

\subsection{One-dimensional cover problem}

\begin{figure}
	\centering
	\includegraphics[width=\linewidth]{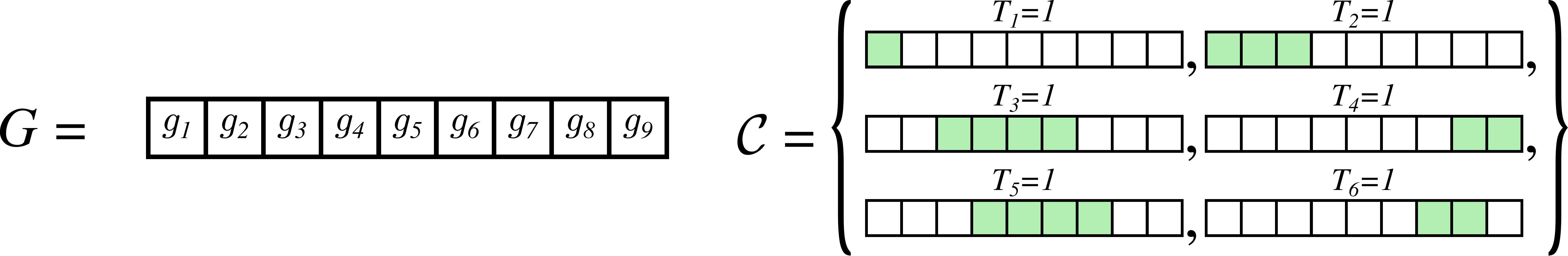}
	\caption{Instance of the one-dimensional cover problem}
	\label{1dcoverproblem}
\end{figure}

For example, let us consider the one-dimensional case of our problem, with $M=1$:
\begin{itemize}
	\item $\grid=\{g_n,\, n\in [0,N[\}$ is a vector with $N$ cells to covers.
	\item a cover $C=\{g_n,\, n\in [n_0,n_1[\}$ is a contiguous subset of $\grid$, uniquely defined by a starting element $n_0$ and an ending element $n_1$ such that $0 \leq n_0 < n_1 \leq N$. Each cover can be represented by a binary vector $$\bc(n) = \begin{cases}
	1 \text{ if } n_0 \leq n < n_1\\  
	0 \text{ otherwise}
	\end{cases}$$ which contains contiguous ``ones", i.e. there is no ``zero" between two ``ones". 
	\item $\covers=\{\cover_1,\dots,\cover_D\}$ is a collection of covers on $\grid$.
\end{itemize}
An example of this problem is presented in Figure \ref{1dcoverproblem}.

For the one-dimensional cover problem, the cover matrix $\bA=\begin{pmatrix}\bc_1 & \cdots & \bc_D\end{pmatrix}$ is an \textit{interval matrix}, i.e. each column $\bc_i$ of $\bA$ has its ones ``consecutively". Interval matrices are known to be unimodular \cite{Nemhauser1999}, and thus totally unimodular, since every submatrix is also an interval matrix.
Every basic optimal solution of the linear relaxation is an integral solution, and a valid solution for the integer program. Solving the linear relaxation of the one-dimensional cover is sufficient to solve the problem itself, making it an ``easy" problem. 

However this not the case for the two-dimensional cover problem. For the problem represented in Figure \ref{setcoverexample} and described by Equation (\ref{setcoverexampleequation}), the linear basic optimal solution is 
$\bx_L=(0\ \tfrac{1}{2}\ \tfrac{1}{2}\ \tfrac{1}{2}\ \tfrac{1}{2}\ \tfrac{1}{2}\ \tfrac{1}{2}\ \tfrac{1}{2})^T$. A possible $\bA_B$ submatrix for this solution is

\begin{figure*}[!h]
\centering
\includegraphics[width=\widthof{$\bA_B=\begin{pmatrix}0&0&0&0&0&0&0&0\end{pmatrix}$}]{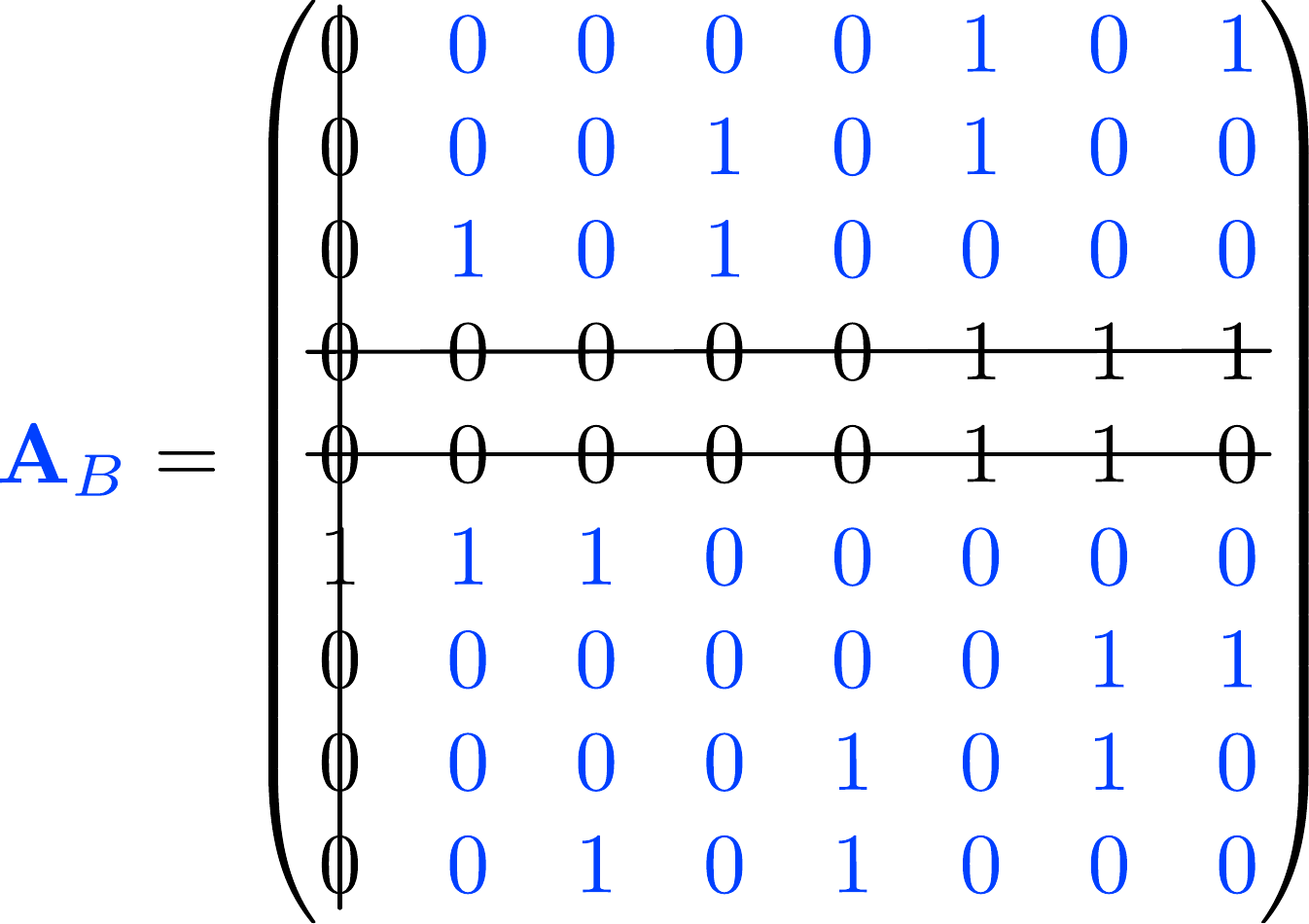}
\end{figure*}

which has a determinant of $-2$, thus explaining the appearance of $\frac{1}{2}$ values in the linear basic solution. This counterexample disproves the total unimodularity of the two-dimensional cover.

\subsection{Integrality Gap}

For the two-dimensional cover problem, the linear optimal solution might not be a valid solution for the integer problem and may have fractional values. The optimal cost of the linear relaxation described by Equation (\ref{setcoverexampleequation}) is $\bT^T\bx_L = \frac{11}{2}$. The optimal cost of the integer program cannot be fractional and is necessarily integer (in this case, it is $6$). The difference between the optimal cost of the integer program and its linear relaxation is called the \textit{integrality gap}.

\subsection{Dynamic programming}

Interestingly, the difficulty gap between the one-dimensional and two-dimensional cover problems can be found by a completely different, algorithmic approach using dynamic programming. Dynamic programming is a method for solving an optimization problem by recursively solving smaller sub-problems. Problems solved by dynamic programming usually possess an \textit{optimal substructure}, which means that an optimal solution can be constructed by combining optimal solutions of its sub-problems. The method is particularly efficient if this substructure can be broken down recursively in a polynomial number of sub-problems.

\begin{figure}
	\centering
	\includegraphics[width=\linewidth]{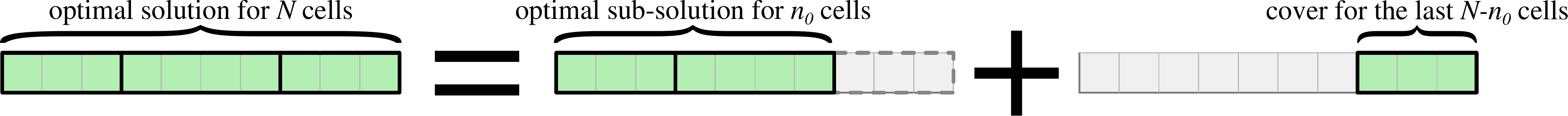}
	\caption{Optimal substructure of the one-dimensional cover problem}
	\label{1doptimalsubstructure}
\end{figure}

This the case for the one-dimensional cover (Fig. \ref{1doptimalsubstructure}). An optimal subcollection $\scanpattern$ covering $\grid$ for this problem is going to be the combination of:
\begin{itemize}
	\item a cover $\bc$ including the last cell, so such that $n_0\leq n < N (=n_1)$. Thus $\bc$ covers the last cell but might also covers some previous cells.
	\item a subcollection covering optimally the first $n_0$ cells (which are not covered by $\bc$)
\end{itemize}

Then it's possible to recursively define the solution covering optimally the first $K$ cells as the union of a one cover including the $K$-th cell and a solution covering optimally the first $k$ cells with $k<K$. This is formalized by the following equation, called the \textit{Bellman recursion}:

$$ \bx_{[1,K]} = \left\{
\begin{array}{cl}
\bold{0} &\text{ if } K=0 \\
\bx_{[1,K-1]} &\text{ if } \bA \bx_{[1,K-1]} \geq \bold{1}_{[1,K]} \\
\argmin \{ \bT^T \bx\ \text{ s.t. }\bx=(\bx_{[1,k]}+\bx_i),\ k<K,\ \bA \bx \geq \bold{1}_{[1,K]}\}&\text{otherwise}
\end{array}\right.
$$

where \begin{itemize}
\item[-] $\bx_{[1,K]}$ is the solution covering optimally at least the first $K$ cells 
\item[-] $\bold{1}_{[1,K]}$ is the vector of length $N$ starting with $K$ ones and ending with $N-K$ zeroes
\item[-] $\bx_i$ is the vector of length $D$ with a one at position $i$ and zeroes elsewhere, representing the addition of the cover $\bc_i$.
\item[-] $k$ is the first non-zero element of the given cover $\bc_i$
\end{itemize}

And the recursion can be described by the following steps

\begin{itemize}
	\item[-] $\bx_{[1,0]}$ is the solution covering zero cells, initialized at $\bold{0}$
	\item[-] If the solution covering $K-1$ cells also covers the $K$-th cell, we keep it.
	\item[-] Otherwise, search the optimal solution covering the first $K$ cells as a combination of: 
	\begin{itemize}
		\item[$\cdot$] a cover $\bc_i$ containing the $K$-th cell
		\item[$\cdot$] the sub-solution $\bx_{[1,k]}$ optimally covering the first $k<K$ cells not covered by $\bc_i$
	\end{itemize}
\end{itemize}

So the dynamic programming method computes optimal sub-solutions for all sub-problems of covering the first $K$ cells for $K \leq N$. So we must solve $N$ sub-problems. For each sub-problem, the solution is computed as a combination of a cover and a smaller substructure optimal solution, so requires $O(D)$ steps to search through all covers. Thus, the algorithmic complexity of the dynamic programming algorithm is $O(ND)$, which is polynomial.

A natural question would be whether this approach can be generalized to the two-dimensional cover problem. Let us consider an optimal solution for the two-dimensional cover problem. It can be viewed as a combination of a rectangular cover $\cover$ including the last bottom-right cell and an optimal cover for the substructure of cell not covered by $\cover$. 

\begin{figure}
	\centering
	\includegraphics[width=\linewidth]{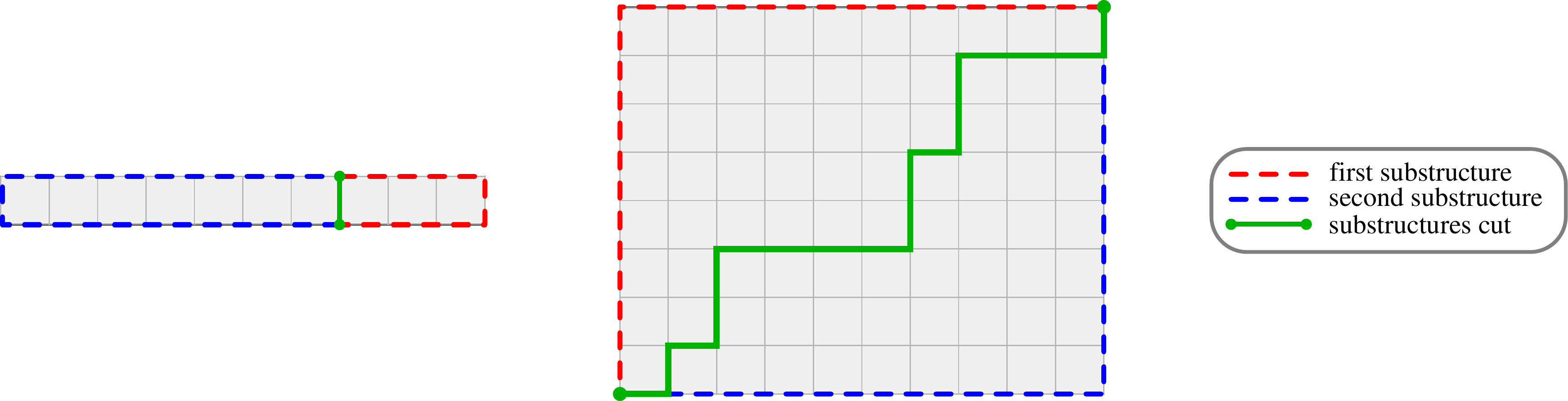}
	\caption{Substructure decomposition of the one-dimensional cover problem (left) and the two-dimensional cover problem (middle)}
	\label{1dvs2doptimalsubstructure}
\end{figure}

So the optimal substructure of the two-dimensional cover problem is the cover sub-problem of a first ``top-left" half of the $M$-by-$N$ grid $\grid$. The number of sub-problems is equal to the number of way of cutting $\grid$ into two substructures: a top-left part and a bottom-right part, as presented in Figure \ref{1dvs2doptimalsubstructure}. Equivalently, this is equal to the number of paths between the top-right corner and the bottom-left corner of $\grid$.

A cut is constituted by $N+M$ edges on the grid, with $M$ vertical edges and $N$ horizontal edges. Any cut can be defined uniquely by choosing the $N$ vertical edges (or equivalently $M$ horizontal edges) among the $N+M$ edges. So the number of possible paths between two opposite corners of $\grid$, and thus the number of cover sub-problems on $\grid$ is $\binom{N+M}{N}=\binom{N+M}{M}$.

Let $K=\min\{N,M\}$, then the number of possible cuts  can be bound below by the following approximation using Stirling's formula
$$\binom{N+M}{N} \geq \binom{2K}{K} \simeq \frac{\sqrt{2\pi2K}(2K)^{2K}}{e^{2K}} \left(\frac{e^{K}}{\sqrt{2\pi K}K^K}\right)^2 = \frac{2^{2K}}{\sqrt{\pi K}}  $$
Thus, the number of sub-problems to solve grows exponentially with the grid size: an increase by $10$ of the grid size increase the number of sub-problems by approximately $2^{2\cdot 10} \approx 10^6$. Even for small values, the number of sub-problems explodes:

\begin{table}[!h]
\centering
\begin{tabular}{|c|c|c|c|c|c|}
	\hline
$N=M$ & 10 & 20 & 30 & 40 & 50 \\
\hline 
$\binom{2N}{N}$ & $\simeq 10^{5}$ & $\simeq 10^{11}$ & $\simeq 10^{17}$ & $\simeq 10^{23}$ & $\simeq 10^{29}$ \\
\hline 
\end{tabular}
\caption{Number of sub-problems}\label{tab}
\end{table}

So while theoretically usable for the two-dimensional cover problem, dynamic programming has an exponential complexity for this problem, making the approach rather unpractical. This hints that the two-dimensional cover problem is computationally harder than the one-dimensional cover problem. The two-dimensional cover problem is in fact NP-difficult to solve \cite{Briheche2018}. To efficiently solve the two-dimensional cover problem, we need to use a more general optimization method.

\section{Branch\&Bound}

Integer programs are generally NP-hard optimization problems: there is currently no know algorithm capable of finding quickly an optimal solution. Informally, exact algorithms ``have to" search through the solution space.

\begin{figure}
	\centering
	\includegraphics[width=\linewidth]{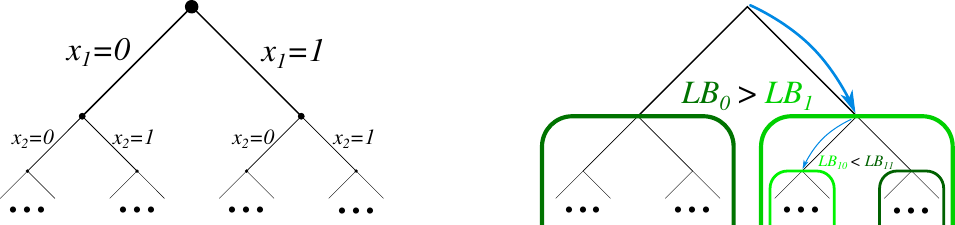}
	\caption{Finite tree of solutions (left) and branch-and-bound method (right)}
	\label{branchbound}
\end{figure}

The space of all possible solutions can be represented as a finite binary tree with depth $p$, each node representing the value choice of an integer variable (Fig. \ref{branchbound}). Each end leaf represents a solution for the integer program. The number of possible solutions is finite, but grows exponentially and is usually huge: in our case, $2^D$ possible solutions.

Exploring the entire tree is computationally unfeasible in reasonable time. However it is possible at each node to estimate a lower bound of the node sub-tree best solution, by solving its linear relaxation. Knowing their lower bound, it is possible to avoid exploring certain subsets. This method is known as the branch-and-bound method \cite{Conforti2014}: 

\begin{itemize}
	\item Branching: Each branch at the current node (with depth $i-1$) correspond to a chosen value, $0$ or $1$, for the next variable $x_i$. In each branch, $x_i$ is no longer a variable but a parameter. The current problem is thus divided into $2$ smaller sub-problems, each considering a different value for $x_i$ and each having one less variable.
	\item Bounding: The current problem is relaxed into a linear program, whose solution is a lower bound of the current problem best solution. Depending on the lower bound value, the node sub-tree will be explored next (if it is the most promising branch), later (if there is a more promising branch), or never (if a better solution has already be found in another branch). 
\end{itemize}

Defining what a promising branch is a difficult question, a lower bound is not necessarily better since deeper nodes may have higher bounds while being closer to optimal solutions. Integer programming solvers usually rely on various heuristics to define the exploration strategy and improve bound estimations.

\subsection{Description}

We present in this section a pseudo-code describing a basic implementation of the branch-and-bound method in Algorithm \ref{branchandbound}. Each node in the tree can be described by the sequence of choices leading to it:
$$\node=(x_1,x_2,\dots,x_n)$$
From a given node, we can compute its children $\node_0 = (x_1,\dots,x_d,0)$ and $\node_1 = (x_1,\dots,x_d,1)$.
At each node $\node$ explored, $(x_1,x_2,\dots,x_n)$ are set, and we solve a linear relaxation of the problem with respect to the variables $(x_{n+1},\dots,x_D)$ using the simplex method, then add $\node$ to the list of nodes to explore. 

The algorithm can be summarized by the following steps:
\begin{itemize}
	\item[\textbf{0.}] \textbf{Initialization}:\\ Initialize the list of node to explore with the root node.
	\item[\textbf{1.}] \textbf{Exploration}:\\
	Pop next node to explore from the list of nodes and compute its linear relaxation.
	\item[\textbf{2.}] \textbf{Bounding}:\\
	If the current node relaxation value is less than the current best solution found, proceed to Step 3, otherwise, drop current node and go back to Step 1.
	\item[\textbf{3.}] \textbf{Update}:\\
	If the current node relaxation is an integral solution, then its an improving solution (note that an end leaf always yield an integral solution). Update best current solution and proceed to Step 1.\\
	Otherwise:
	\item[\textbf{4}.] \textbf{Branching}:\\ Compute the current node children. For each child, check if the descendants contains a valid solution (this can be done by summing covers already used by the parent, the cover of the child node if used, and covers available to the descendants). If the child node is valid, add it to the list of node to explore. Proceed to Step 1.
\end{itemize}

\begin{algorithm}[!h]
	\caption{Branch-and-bound}\label{branchandbound}
	\begin{algorithmic}
		\State \textcolor{gray}{\% \Call{lp\_solve}{} is the relaxation subroutine called during branching}
		\Function{lp\_solve}{$\node$}
		\State $(x_1,\dots,x_{d-1}) := \node$ \Comment{At node $\node$, the first $d-1$ variables are set}
		\State $(x_d,\dots,x_D) := \argmin\{\sum_{j=d}^{D}T_j x_j \text{, s. t. } \bA \cdot \bx \geq \bold{1}\}$ \Comment{Optimization of non-set variables}
		\State \Return $\bx_L := (x_1,\dots,x_d,x_{d+1},\dots,x_D)$
		\EndFunction

		\\ \State \textcolor{gray}{\% Initialization}
		\State $\node_{root}=()$ 
		\State $\nodes := \{\node_{root}\}$ \Comment{Start with root node}
		\State $\bx_{best} := \bx_F = (1\cdots1)$ \Comment{Best solution found so far (by default, $\bx_F$ is a valid solution)}

		\\ \State \textcolor{gray}{\% Exploration}
		\While{$\nodes$ is not empty}

		\State $\node:=\operatorname{pop}(\nodes)$ \Comment{Take next node in $\nodes$}
		\State $\bx_L := \Call{lp\_solve}{\node}$  \Comment{Solve node relaxation}

		\\ \State \textcolor{gray}{\% Bounding}
		\If{ $\bT^T\cdot \bx_L < \bT^T \cdot \bx_{best}$ } \Comment{Explore node $\node$ only if it can improve best solution}
			\\ \State \textcolor{gray}{\% Update}
			\If { $\bx_L \in \{0,1\}^D$ } \Comment{Check if $\bx_L$ is an integral solution}
				\State $\bx_{best} := \bx_L$
			\Else
				\State $(x_1,\dots,x_d) := \node$\\
				
				\State \textcolor{gray}{\% Branching}
				\For {$x \in \{0,1\}$}  \Comment{Compute children of node $\node$}
					\State $\node_c := (x_1,\dots,x_d,x)$					
					\If {$\sum_{j=1}^{d}x_j\bc_j+x\bc_{d+1}+\sum_{j=d+2}^{D}\bc_j\ \geq\ \bold{1}$} \Comment{Check child feasibility}
						\State $\nodes := \nodes \cup \{\node_c\}$ \Comment{Add child to the candidate list}
					\EndIf
				\EndFor
			\EndIf

		\EndIf
		\EndWhile
		\State \textbf{return} $\bx_{best}$
	\end{algorithmic}
\end{algorithm}

This very generic description is just a presentation of the general idea of the method. Efficient implementations of the branch-and-method usually combined several techniques such as cutting planes, diving heuristics and local branching to improve bounds estimation and speed.

%
%

\subsection{Application example}

In this section, we apply and describe the behavior of the branch-and-bound method on the example represented in Figure \ref{setcoverexample} and described in Equation \ref{setcoverexampleequation}.

\begin{itemize}
	\item $\nodes=\{\}$, $\bx_{best}=(1\ 1\ 1\ 1\ 1\ 1\ 1\ 1)$, $f_{best}=\bT^T\cdot\bx_{best}=13$ :\\
	Solving the root relaxation yields the linear solution $(0\ \tfrac{1}{2}\ \tfrac{1}{2}\ \tfrac{1}{2}\ \tfrac{1}{2}\ \tfrac{1}{2}\ \tfrac{1}{2}\ \tfrac{1}{2})$ with cost $\tfrac{11}{2}\leq 13$. Root node children $(0)$ and $(1)$ are feasible, and thus added to the exploration list $\nodes:=\{(0),(1)\}$
	
	\item $\nodes=\{(0),(1)\}$, $\bx_{best}=(1\ 1\ 1\ 1\ 1\ 1\ 1\ 1)$, $f_{best}=\bT^T\cdot\bx_{best}=13$ :\\
	Relaxation of $(0)$ yields the same linear solution $(0\ \tfrac{1}{2}\ \tfrac{1}{2}\ \tfrac{1}{2}\ \tfrac{1}{2}\ \tfrac{1}{2}\ \tfrac{1}{2}\ \tfrac{1}{2})$ with cost $\tfrac{11}{2}$. We add the children $(0,0)$ and $(0,1)$ to the exploration list $\nodes:=\{(1),(0,0),(0,1)\}$
	
	\item $\nodes=\{(1),(0,0),(0,1)\}$, $\bx_{best}=(1\ 1\ 1\ 1\ 1\ 1\ 1\ 1)$, $f_{best}=\bT^T\cdot\bx_{best}=13$ :\\
	Relaxation of $(1)$ yields the linear optimal solution $\bx_L=(1\ 0\ 1\ 1\ \tfrac{1}{2}\ \tfrac{1}{2}\ \tfrac{1}{2}\ \tfrac{1}{2})$ with cost $\tfrac{15}{2}<13$. We add the children $(0,0)$ and $(0,1)$ to the exploration list $\nodes:=\{(1,0),(1,1)\}$

	\item $\nodes=\{(0,0),(0,1),(1,0),(1,1)\}$, $\bx_{best}=(1\ 1\ 1\ 1\ 1\ 1\ 1\ 1)$, $f_{best}=\bT^T\cdot\bx_{best}=13$ :\\
	Relaxation of $(0,0)$ yields the linear optimal solution $\bx_L=(0\ 0\ 1\ 1\ 0\ 0\ 1\ 1)$ with cost $6<13$. $\bx_L$ is an integral solution, thus we update the best current solution $\bx_{best}:=\bx_L;\ f_{best}:=6$.	
\end{itemize}

At this point, we can already deduce than we have found an integer optimal solution. The root relaxation has linear optimal cost $\tfrac{11}{2}$. Since any integer solution is a valid linear solution, it has an integer cost greater than the linear optimal cost $\tfrac{11}{2}$, so greater than $6$. This suffices to prove the optimality of $\bx_{best}=(0\ 0\ 1\ 1\ 0\ 0\ 1\ 1)$ for the integer program described by Equations (\ref{integerprogram},\ref{setcoverexampleequation}).

\subsection{Multiple solutions enumeration}

While we could terminate the exploration once we have found an optimal solution, we also have the possibility to pursue the exploration in order to found alternative optimal solutions. 

In engineering applications, multiple solutions are a desirable feature for engineers and operators who can select a solution among multiple candidates based on their expertise. This choice in turn can be analyzed to define preferences, to add secondary selection criterion to the method or even refined the model into a multi-objective optimization problem. 

Multiple solutions enumeration can be done by slightly modifying steps \textbf{2.} and \textbf{3.} of the branch-and-bound method:

\begin{itemize}
	\item[\textbf{2.}] \textbf{Bounding}:\\
	If the current node relaxation value is less than or equal to the current best solution found, proceed to Step 3, otherwise, drop current node and go back to Step 1.
	\item[\textbf{3.}] \textbf{Update and Enumerate}:\\
	If the current node relaxation is an integral solution, then its an improving solution. If it is strictly better than the current solution, empty the set of best solutions and update best current solution. Otherwise, update the set of best solutions. Proceed to Step 4 (as there could be other optimal solutions among the children of the current node).\\
\end{itemize}

This result in modifications to Algorithm \ref{branchandbound} pseudo-code as described in Algorithm \ref{enumerationbranchandbound}.
\begin{algorithm}[!h]
	\caption{Multiple solutions enumeration branch-and-bound}\label{enumerationbranchandbound}
	\begin{algorithmic}
		\State \textcolor{gray}{\% Initialization}
		\State ...
		\State $\bx_{best} := \bx_F = (1\cdots1)$ \Comment{Best solution found so far (by default, $\bx_F$ is a valid solution)}
		\State $\mathcal{X}_{best} := \{\bx_F\}$ \Comment{Set of best solutions found so far}

		\\ \State \textcolor{gray}{\% Exploration}
		\While{$\nodes$ is not empty}
		\State ...
		
		\\ \State \textcolor{gray}{\% Bounding}
		\If{ $\bT^T\cdot \bx_L \leq \bT^T \cdot \bx_{best}$ } \Comment{Explore $\node$ if its relaxation is at least as good as $\bx_{best}$}
		\\ \State \textcolor{gray}{\% Update and Enumerate}
		\If { $\bx_L \in \{0,1\}^D$ } \Comment{Check if $\bx_L$ is an integral solution}
		\If{ $\bT^T\cdot \bx_L < \bT^T \cdot \bx_{best}$ }
		\State $\bx_{best} := \bx_L$
		\State $\mathcal{X}_{best} := \{\bx_L\}$
		\Else
		\State $\mathcal{X}_{best} := \mathcal{X}_{best} \cup \{\bx_L\}$
		\EndIf
		\EndIf \\
		\State \textcolor{gray}{\% Branching}
		\For {$x \in \{0,1\}$}
		\State ...
		\EndFor			
		\EndIf
		\EndWhile
		\State \textbf{return} $\mathcal{X}_{best}$
	\end{algorithmic}
\end{algorithm}

If we pursue the method application to the numerical example previously described:
\begin{itemize}
	\item $\nodes=\{(0,0),(0,1),(1,0),(1,1)\}$, $\bx_{best}=(1\ 1\ 1\ 1\ 1\ 1\ 1\ 1)$, $f_{best}=\bT^T\cdot\bx_{best}=13$ :\\
	Relaxation of $(0,0)$ yields the linear optimal solution $\bx_L=(0\ 0\ 1\ 1\ 0\ 0\ 1\ 1)$ with cost $6\leq13$. $\bx_L$ is an integral solution, thus we update the best current solution $\bx_{best}:=\bx_L;\ f_{best}:=6$.\\
	We add the children $(0,0,0)$ and $(0,0,1)$ to the exploration list $\nodes$.

	\item $\nodes=\{(0,1),(1,0),(1,1),(0,0,0),(0,0,1)\}$, $\bx_{best}=(0\ 0\ 1\ 1\ 0\ 0\ 1\ 1)$, $f_{best}=\bT^T\cdot\bx_{best}=6$ :\\
	Relaxation of $(0,1)$ yields the linear optimal solution $\bx_1=(0\ 1\ 1\ 0\ 0\ 1\ 1\ 0)$ with cost $6\leq6$. $\bx_1$ is an integral solution, thus added to $\mathcal{X}_{best}:=\{\bx_{best},\bx_1\}$.
	We add the children $(0,1,0)$ and $(0,1,1)$ to the exploration list $\nodes$.
		
	\item $\nodes=\{(1,0),(1,1),(0,0,0),\dots\}$, $\bx_{best}=(0\ 0\ 1\ 1\ 0\ 0\ 1\ 1)$, $f_{best}=\bT^T\cdot\bx_{best}=6$ :\\
	Relaxation of $(1,0)$ yields the linear optimal solution $\bx_L=(1\ 0\ 1\ 1\ \tfrac{1}{2}\ \tfrac{1}{2}\ \tfrac{1}{2}\ \tfrac{1}{2})$ with cost $\tfrac{15}{2}>6$. We drop node $(1,0)$ and proceed with the next node.

	\item $\nodes=\{(1,1),(0,0,0),\dots\}$, $\bx_{best}=(0\ 0\ 1\ 1\ 0\ 0\ 1\ 1)$, $f_{best}=\bT^T\cdot\bx_{best}=6$ :\\
	Relaxation of $(1,1)$ yields the linear optimal solution $\bx_L=(1\ 1\ 1\ 0\ 0\ 0\ 1\ 1)$ with cost $8>6$. We drop node $(1,1)$ and proceed with the next node.
\end{itemize}

\begin{figure}
	\centering
	\includegraphics[width=0.6\linewidth]{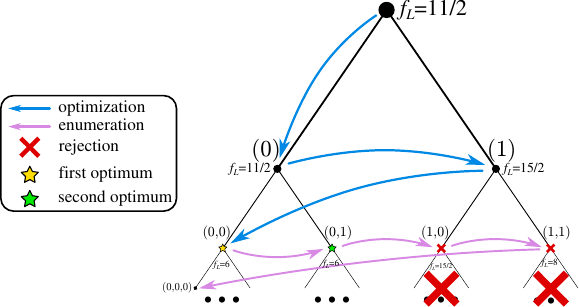}
	\caption{Graphical representation of the branch-and-bound application example}
	\label{branchandboundapplication}
\end{figure}

This numerical example is graphically represented in Figure \ref{branchandboundapplication}, with the optimization phase, the enumeration phase and some nodes rejection.

\subsection{Just-in-time criteria}

One of the most interesting features of the branch-and-bound method from an operational point of view is the possibility to use a ``just-in-time" criteria. For example, a radar system with an embedded computer must optimize its cover just before a mission start. However, it only has five minutes to perform the optimization. A ``just-in-time" is a time limit condition that would ensure that even if the optimum has not been reached, the algorithm will return the best solution it found in the available lapse of time.
Another strength of the method is the fact that linear relaxation provides a lower bound of the optimal solution value: 
$$B_{\nodes}=\min\{\bT^T\cdot \bx_L\ :\ \bx_L=\text{LP\_SOLVE}(\node), \node \in \nodes \}$$
thus during the computation of the method, we always have an interval of confidence for the optimal solution value, above the lower bound but below the current best value:
$$B_\nodes \leq \bT\cdot \bx_{opt} \leq \bT\cdot \bx_{best}$$
Knowing the lower bound, we can compute the \textit{(worst-case) relative optimality gap} as:
$$\Delta_{opt} = \frac{\bT\cdot \bx_{best} - B_\nodes}{B_\nodes}$$
which give as a percentage the best gain we can hope from the optimal solution relatively to the current best solution. The pseudo-code modifications required to account a time limit and provided the current lower bound are described in Algorithm \ref{jitbranchandbound}.

\begin{algorithm}[!h]
	\caption{Just-in-time branch-and-bound}\label{jitbranchandbound}
	\begin{algorithmic}
		\\ \State \textcolor{gray}{\% Exploration}
		\State current\_time := time() \Comment{Get current time}
		\While{$\nodes$ is not empty AND current\_time $\leq$ time\_limit }
		\State ...
		
		\EndWhile
		\State \textbf{return} $\mathcal{X}_{best}$, $B_\nodes$
	\end{algorithmic}
\end{algorithm}

In practice, if the algorithm has a broad choice of available covers, it will find very quickly a good quality solution. Typically within $\leq 10\%$ of relative optimality gap. However closing those last percents to reach the optimal solution can be difficult. Because the decision space is often huge, the algorithm spends a long time crossing out possibilities. In some case even, the algorithm finds quickly the optimal solution, and spends a long time proving its optimality. 

\section{Application to Radar Engineering}

In this section, we give a study case example of radar search pattern optimization and its simulation results. We present first a quick informal description with intuitive and quantitative insight on our mathematical model of the radar system.

\subsection{Radar model}

An active radar is a system capable of detecting distant metallic objects, by sending electromagnetic waves and listening to reflected echos. To perform detection in a given azimuth-elevation direction $(az,el)$, the radar antenna is electronically controlled to focus power in direction $(az,el)$, maximizing the radiation pattern in that direction. A signal containing a series of impulses is then sent through the radar. Upon reception, the reflected signal is filtered to detect echoes. A longer signal is more energetic and easier to filter out. The energy received by the radar from a target at distance $R$ is
\begin{equation}
	E_r = K\ \frac{g^2\ T}{R^4}
	\label{radarequation}
\end{equation}
where $K$ is a constant accounting for the radar emitting power, internal losses, target reflectability, etc., $g$ is the antenna radiation pattern in direction $(az,el)$, and $T$ is the signal duration.

Equation (\ref{radarequation}) is a simpler version of the \textit{radar equation} \cite{Skolnik2008}, a fundamental concept in radar theory. It formalize the intuitive idea that the reflected energy increases with antenna directivity and signal duration, but decreases with the target distance. The radar has a certain detection threshold $E_t$, and detects a target only if its reflected energy is above this threshold.

For a given radar system, we have a set of feasible rectangular radiation patterns for the antenna and a set of available signals. The combination of a radiation pattern and a time signal is called a \textit{dwell}.

\subsection{Simulation parameters}

\begin{figure}
	\centering
	\includegraphics[width=\linewidth]{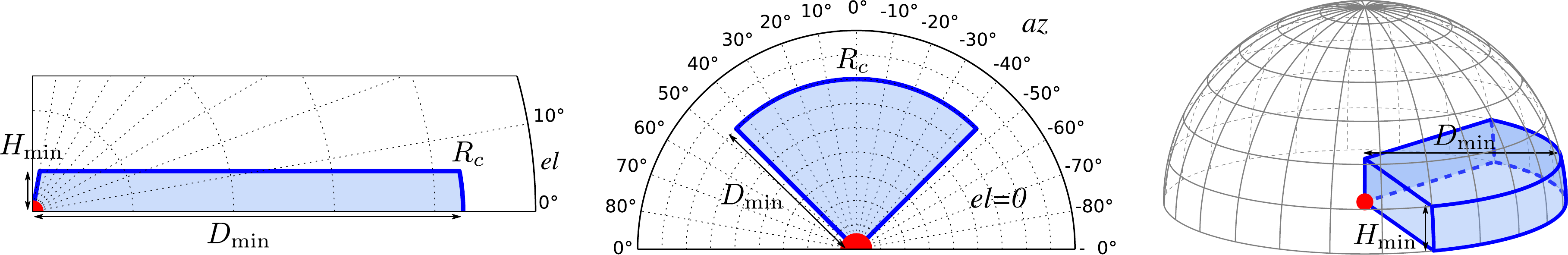}
	\caption{Desired detection range: elevation cut (left), azimuth cut (center) and 3D view (right)}
	\label{desiredrange}
\end{figure}

The desired detection range is defined by a minimum distance $D_{\min}$and a minimum altitude $H_{\min}$. We want to detect targets within that range, closer than $D_{\min}$ and below altitude $H_{\min}$, so the desired detection range is defined as:
$$\range_c(az,el)=
\left\{\begin{array}{cl}
D_{\min} & \text{ if } el \leq \operatorname{asin}\left(\frac{H_{\min}}{D_{\min}}\right)\\
\frac{H_{\min}}{\sin(el)} & \text{ otherwise }
\end{array}\right.$$
Informally, the volume defined by the detection range resembles a sliced cylinder (Fig. \ref{desiredrange}). The radar can use two different type of signals: a short signal and a long signal. 

In our simulation, the detection grid $G$ is a 20$\times$20 lattice with 326 valid cells. We computed 866 feasible dwells in our study case, with 815 dwells using a short time signal with duration $T_s$ and 51 dwells using a long time signal with duration  $T_l$. We compute the cost vector $\bT$ with size 866 associating each dwell to its signal time duration.

We want to find a optimal radar search pattern, i.e. a sub-collection of dwells among the 866 available dwells covering all 326 valid detection cells with minimal total time-budget.  For each of the 866 dwells, we use Equation (\ref{radarequation}) to compute the dwell detection cover on the 326 cells. From the detection covers, we can compute the cover matrix $\bA$ with shape 326 $\times$ 866.

Having computed $\bT$ and $\bA$, we can use the branch-and-bound method described previously to search an optimal radar search pattern. The corresponding integer program has 866 variables and 326 detection constraints.

The optimization is done through the CPLEX solver \cite{CPLEX}, which implements an improved version of the branch-and-bound. The total time required to find the solution is 5 seconds on an i7-3770@3.4GHz processor.

\begin{figure}
	\centering
	\includegraphics[width=\linewidth]{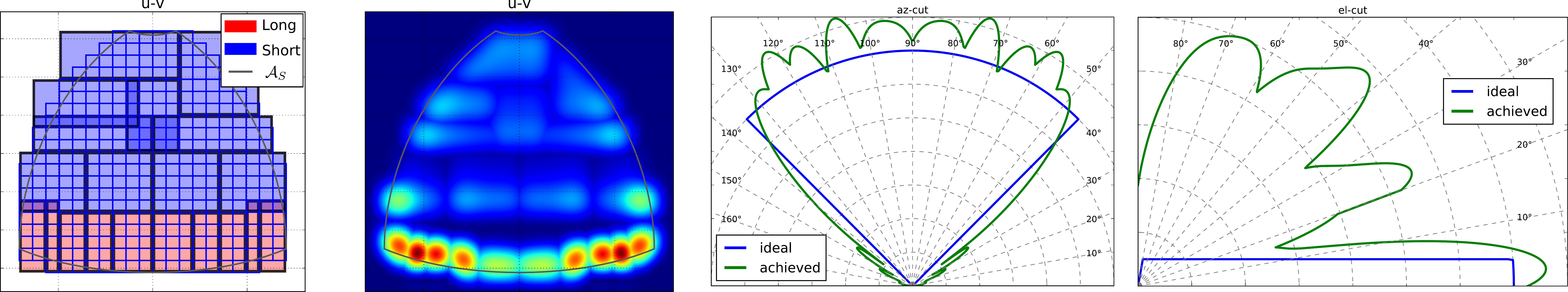}
	\caption{From left to right: computed radar search pattern, emission pattern and detection range}
	\label{solution}
\end{figure}

\subsection{Optimal solution}

The returned optimal solution is shown in Figure (\ref{solution}): the left sub-figure shows the discrete covers of the 20 dwells used in the pattern. 10 dwells use a short signal (represented in blue) and cover high elevation, and 10 dwells use a long signal (represented in red) and cover low elevations. 

This result is explained by the fact that a radar must usually achieve high detection range near the horizon (where targets are located) and low detection range at high elevation (since most aircrafts have limited flying altitude). It makes sense to use longer, and thus ``more energetic" dwells at low elevations than at high elevations.

\subsection{Enumeration}

As we have seen before, there may be multiple optimal solutions. In this simulation we managed to found 3500 different optimal solutions in 5 minutes. However, this search is unlikely to be exhaustive: due to the wide choice of possible dwells, there is often an extremely high number of possible alternative optimal solutions. Finding all solutions is unfeasible in practice.

However optimal solutions share certain characteristics: all solutions have 10 shorts signal dwells at high elevation and 10 long signal dwells at low elevation. The long-signal dwells (in red) are mostly the same for all optimal solutions found, and form an \textit{optimality invariant}. The short-signal dwells (in blue) are however different for each solution.

Intuitively, the low-elevation area is more ``energetically demanding"; thus low-elevation detection constraints are the ``hardest constraints" of the problem, and do not leave a lot of choice for covering the low-elevation area. High-elevation detection constraints are in comparison ``easier" and can be validated by different covers.

\begin{figure}[b]
	\centering
	\includegraphics[width=0.5\linewidth]{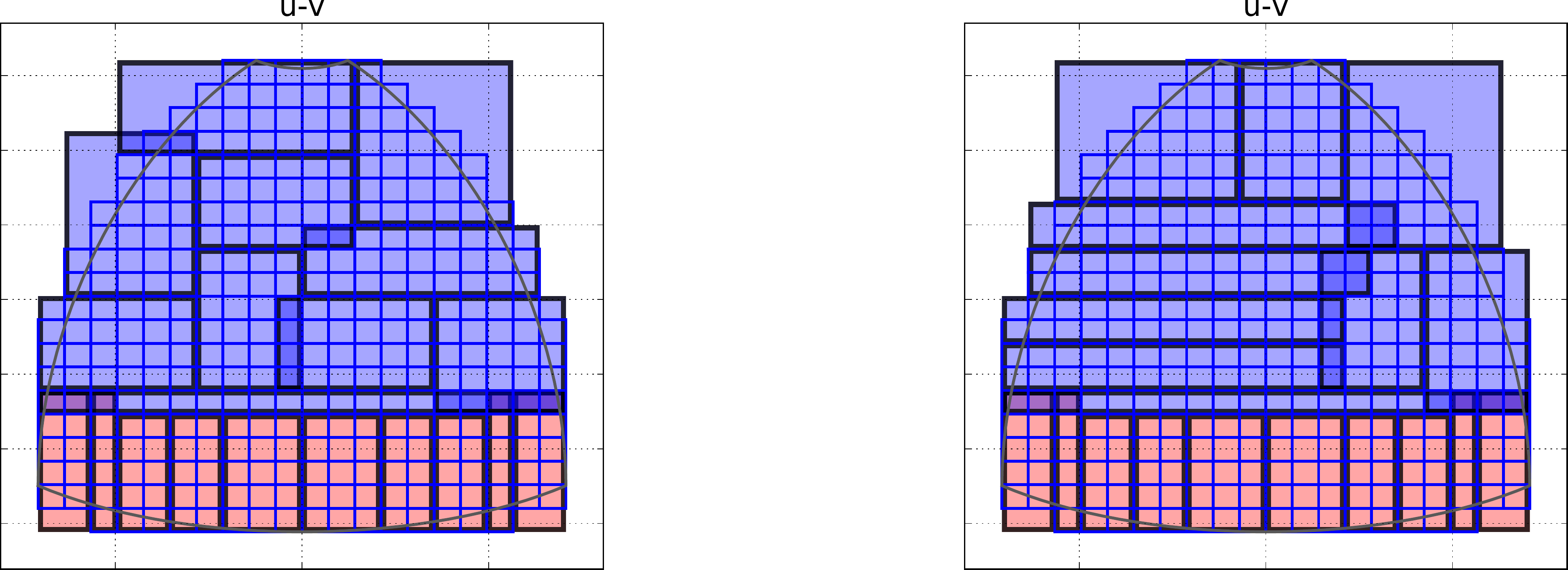}
	\caption{Two other possible optimal solutions}
	\label{alternatesolutions}
\end{figure}

\section{Conclusion}

The branch-and-bound method is a practical and powerful technique. It can be used as an exact algorithm if exploration is pushed to its completion, when there is no more branch left with a potentially better solution. It can also be used as an heuristic, with stopping criterion based on an time limit or a optimality gap treshold. This is especially useful in operational situations with broad choices, when finding a good solution is easy, but proving optimality is difficult.

The method is very generic, and can be used to solve a lot of different combinatorial problems. Many of those problems have evident practical values and important applications in various industries, such as the set cover problem in radar applications. The versatility of the method and its various ``flavors" can be used for different purposes: enumeration permits analysis of the radar ``possibilities" during conception, while just-in-time criteria improves resources management in operational situations. Branch-and-bound is an extremely efficient tool for a broad variety of engineering applications.


\section*{Acknowledgments}
This work is partly supported by a DGA-MRIS scholarship.

\bibliographystyle{ietConference}
\bibliography{biblio}

\end{document}